\documentclass[12pt,reqno]{amsart}
\usepackage{amsthm,amsfonts,amssymb,amscd}
\newtheorem{theorem}[subsection]{Theorem}
\newtheorem{proposition}[subsection]{Proposition}
\newtheorem{lemma}[subsection]{Lemma}
\newtheorem{corollary}[subsection]{Corollary}
\theoremstyle{definition}
\newtheorem{definition}[subsection]{Definition}
\newtheorem{problem}[subsection]{Problem}
\newtheorem{example}[subsection]{Example}
\theoremstyle{remark}
\newtheorem{remark}[subsection]{Remark}
\newcommand{\mt}[1]{\operatorname{#1}}
\newcommand{\codim}{\operatorname{codim}}
\newcommand{\NE}{\operatorname{NE}}
\newcommand{\Diff}{\operatorname{Diff}}
\newcommand{\Supp}{\operatorname{Supp}}

\newcommand{\Sing}{\operatorname{Sing}}
\newcommand{\bir}{\dasharrow}

\newcommand{\down}[1]{\lfloor #1\rfloor}
\newcommand{\fr}[1]{\{ #1\}}
\renewcommand{\AA}{{\mathbb A}}
\newcommand{\CC}{{\mathbb C}}

\newcommand{\ZZ}{{\mathbb Z}}
\newcommand{\QQ}{{\mathbb Q}}
\newcommand{\PP}{{\mathbb P}}
\newcommand{\NN}{{\mathbb N}}

\newcommand{\OOO}{{\mathcal O}}

\renewcommand\labelenumi{(\roman{enumi})}

\email{prokhoro@mech.math.msu.su\qquad \\
prokhoro@math.pvt.msu.su\\}
\begin{document}
\author{Yu.~G.~Prokhorov}
\address{Department of Mathematics
(Algebra Section), Moscow State University, 117234
Moscow, Russia}
\title{Blow-ups of canonical singularities}
\maketitle
\section*{Introduction}
In this paper we apply Shokurov's inductive method
to study terminal and canonical singularities.
As an easy consequence of the Log Minimal Model Program
we show that for any three-dimensional log terminal singularity
$(X\ni P)$ there exist some special, so called,
plt blow-ups $Y\to X$ (Proposition~\ref{plt}).
We discuss properties of them and construct some examples.
Roughly speaking a plt blow-up $Y\to X$ is a blow-up with an
irreducible exceptional divisor $S$ such that
$S$ is normal and the pair $(S,\Diff_S(0))$ is log terminal,
where $\Diff_S(0)$ is the different \cite[\S 3]{Sh}.
These blow-ups are very useful for inductive approach to
the study of singularities and, more general, extremal contractions
\cite{Pro2}.
We also obtain necessary conditions for
a log surface to be an exceptional divisor of
a plt blow-up of a terminal singularity (Corollary~\ref{m1},
Proposition~\ref{cond2}).
Unfortunately a plt blow-up of a terminal
singularity is never unique (see \ref{term-non-un}).
However there are canonical singularities for which it is unique.
These singularities are called weakly exceptional
(see Definition~\ref{def-weakly-exc}) and have the
most interesting and complicated structure (cf. \cite{MP}, \cite{MP1}).
We obtain the criterion for a singularity to be
 weakly exceptional, in terms of the exceptional divisor of
some plt blow-up (Theorem~\ref{weakly-exc}) and construct few examples.

\par {\footnotesize
\textbf{Acknowledgements.}
I am grateful to V.~V.~Shokurov for useful discussions
concerning complements and exceptionality.
I wish to thank  also D.~Markushevich for
many helpful conversations during my stay
at the University of Lille 1 in the spring of 1998.
This work was partially supported by a grant
96-01-00820 from the Russian Foundation of Basic Research
and a grant PECO-CEI from the Ministry of Higher Education
of France. }

\section{Preliminary results}
All varieties are algebraic and are assumed to be defined over
 $\CC$, the complex number field.
We follow the terminology and notation of the log Minimal Model
Program \cite{KMM}, \cite{Sh}, \cite{Ut} (see also \cite{KoP}
for an introduction to singularities of
pairs). We will use the following abbreviations:
\par
\begin{tabular}{ll}
lc& log canonical \cite{KMM}\\
klt& Kawamata log terminal \cite{Sh}, \cite{Ut}, \cite{KoP} \\
plt& purely log terminal \cite{Sh}, \cite{Ut}, \cite{KoP}\\
dlt& divisorial log terminal \cite{Sh}, \cite{Ut}\\
$LCS(X,D)$&the locus of lc singularities of $(X,D)$ \cite[3.14]{Sh},
\cite{K1}\\
\end{tabular}
\par
A \textit{contraction} (or {\it extraction}, if we start with
$X$ instead of $Y$) is a projective morphism of
normal varieties $f\colon Y\to X$ such that $f_*\OOO_Y=\OOO_X$.
A \textit{blow-up} is a birational extraction.
Usually we write $f\colon (Y,S)\to X$
to specify that $S$ is the (reduced) exceptional divisor of $f$.
Two blow-ups $f\colon Y\to X$
and $f'\colon Y'\to X$ are said to be \textit{birationally
equivalent} if the map $(f^{-1}\circ f')\colon Y'\bir Y$
is an isomorphism in codimension one. We will work in the
category of algebraic varieties, hovewer many  results
can be easily modified to the category of analytic spaces.

\begin{theorem}[{\cite[3.1]{Sh2}}, {\cite[17.10]{Ut}},
cf. {\cite[5.19]{Sh}}]
\label{min-lt}
Let $X$ be a normal variety of dimension $\le 3$ and
let $D$ be a boundary on it such that $(X,D)$ is lc.
Then there exists a blow-up $f\colon Y\to X$ such that
\begin{enumerate}
\item
every exceptional divisor $E_i$ of $f$ has the
discrepancy $a(E_i,D)=-1$ (and therefore one can write
$f^*(K_X+D)=K_Y+D_Y+\sum E_i$, where $D_Y$ is the proper
transform of $D$);
\item
$(Y, D_Y+\sum E_i)$ is dlt;
\item
$Y$ is $\QQ$-factorial.
\end{enumerate}
\end{theorem}
We call such a blow-up \textit{minimal log terminal} blow-up
of $(X,D)$.

\begin{definition}[{\cite[\S 3]{Sh}}, {\cite[Ch. 16]{Ut}}]
Let $X$ be a normal variety and let $S+B$ be a boundary
on $X$, where $S=\down{S+B}\ne \emptyset$ and $B=\fr{S+B}$.
Assume that $K_X+S+B$ is lc in codimension two.
Then the \textit{different} of $B$ on $S$ is defined by
$$
K_S+\Diff_S(B)\equiv (K_X+S+B)|_S
$$
\end{definition}

\begin{theorem}[{\cite[3.9]{Sh}}, {\cite[16.6]{Ut}}]
\label{coeff}
Let $(X,S)$ be a plt pair, where $S$ is reduced.
Then $\Diff_S(0)=\sum (1-1/m_i)R_i$,
where $R_i$ are irreducible components of $\Sing(X)$
contained on $S$ and having codimension one in $S$.
The numbers $m_i\in\NN$ can be characterized by one of the
following properties:
\begin{enumerate}
\item
$m_i$ is the index of $K_X+S$ at the general point of $R_i$;
\item
$m_i$ is the index of $S$ at the general point of $R_i$.
\end{enumerate}
\end{theorem}
We say that coefficients of a $\QQ$-divisor $D$
are \textit{standard} if they have the form $(1-1/m_i)$,
where $m_i\in\NN\cup\{\infty\}$.

\begin{theorem}[Inversion of Adjunction
{\cite[3.3]{Sh}}, {\cite[17.6]{Ut}}]
\label{Inv-Adj}
Let $X$ be a normal variety and let $D$
be a boundary on it. Write $D=S+B$, where $S=\down{D}$
and $B=\fr{D}$. Assume that $K_X+S+B$ is $\QQ$-Cartier.
Then $(X,S+B)$ is plt near $S$ iff $S$ is normal and
$(S,\Diff_S(B))$ is klt.
\end{theorem}

\begin{definition}[{\cite[5.1]{Sh}}]
\label{def-compl}
Let $X$ be a normal variety and let $D=S+B$ be a
boundary, where $S=\down{D}$ is the integer part and $B=\fr{D}$
is the fractional part. Then one
says that $K_X+D$ is \textit{$n$-complementary}, if there exists a
$\QQ$-divisor $D^{+}$, such that the following conditions are verified:
\begin{enumerate}
\item
$n(D^{+}+K_X)\sim 0$ in particular $nD^{+}$ has integer coefficients;
\item
$nD^{+}\ge nS+\down{(n+1)B} \geq 0$;
\item
 $K_X+D^{+}$ is lc.
\end{enumerate}
The log canonical divisor $K_X+D^{+}$ is called an
\textit{$n$-complement} of $K_X+D$.
\end{definition}
Note that in general it is not true that $D^+\ge D$.
This however is true if the coefficients of $D$ are standard
\cite[2.7]{Sh1}.

\begin{proposition}[{\cite[5.4]{Sh}}, {\cite[4.4]{Sh1}}]
\label{podnyal}
Let $f\colon Y\to X$ be a blow-up and
let $D$ be a boundary on $Y$. Then
\begin{enumerate}
\item
$K_Y+D$ is $n$-complementary $\Longrightarrow$
$K_X+f(D)$ is $n$-complementary.
\item
Assume additionally that $K_Y+D$ is $f$-nef and
$f(D)$ has only standard coefficients.
Then \par
$K_X+f(D)$ is $n$-complementary $\Longrightarrow$
$K_Y+D$ is $n$-complementary.
\end{enumerate}
\end{proposition}
Combining \cite[19.6]{Ut} with \ref{podnyal} we get the following.

\begin{proposition}[{\cite[2.1]{Pro_conic}}, cf. {\cite[19.6]{Ut}}]
\label{prodolj}
Let $(X,S)$ be a purely log terminal pair with reduced
$S\ne 0$ and let $f\colon X\to Z$ be a projective morphism
such that $-(K_X+S)$ is
$f$-nef and $f$-big. Given an $n$-complement $K_S+\Diff_S(0)^+$ of
$K_S+\Diff_S(0)$, then in a neighborhood of any fiber of
$f$ meeting $S$ there exists
an $n$-complement $K_X+S+B$ of $K_X+S$ such that
$\Diff_S(0)^+=\Diff_S(B)$.
\end{proposition}
$n$-complements with $n\in\{1, 2, 3, 4, 6\}$
are called \textit{regular}.

\begin{proposition}[{\cite{Sh1}}]
\label{DelPezzo}
Let $(S,\Theta)$ be a projective log surface
such that $(S,\Theta)$ is lc but is not klt and $-(K_S+\Theta)$
is nef and big (i.~e. $(S,\Theta)$ is a weak log Del Pezzo surface,
see~\cite{Sh1}). Then $K_S+\Theta$ has a regular non-klt complement
$K_S+\Theta^+$. Moreover,
we can take  $\Theta^+$ so that $a(E,\Theta)=-1$
implies   $a(E,\Theta^+)=-1$   for every (not necessary
exceptional divisor $E$.
\end{proposition}
\begin{proof}
Let $f\colon\widetilde{S}\to S$ be a log resolution of $(S,\Theta)$.
Consider the crepant pull back $K_{\widetilde{S}}+E+D=
f^*(K_S+\Theta)$, where $E$ is a reduced divisor, $D$
is an effective $\QQ$-divisor and $\down{D}\le 0$.
By \cite[5.2]{Sh} $K_E+\Diff_E(D)$ is
$1$, $2$, $3$, $4$ or $6$-complementary
and by \cite[19.6]{Ut} this complement can be extended on
$\widetilde{S}$. Finally, we can push down it on $S$
by \ref{podnyal}.
\end{proof}

The following proposition can be proved in the same style.

\begin{proposition}[{\cite{Sh1}, \cite{Pro}}]
\label{P-bundle}
Let $f\colon S\to Z$ be a contraction from
a surface $S$ to a curve $Z$ and let $\Theta$
be  a boundary on $S$. Assume that
 $(S,\Theta)$ is lc and $-(K_S+\Theta)$
is $f$-ample. Then
near a fiber $f^{-1}(o)$, $o\in Z$ the log divisor
$K_S+\Theta$ has a regular non-klt complement $K_S+\Theta^+$.
Moreover, if  $(S,\Theta)$ is  not klt, then
we can take  $\Theta^+$ so that $a(E,\Theta)=-1$
implies   $a(E,\Theta^+)=-1$   for every (not necessary
exceptional divisor $E$.
\end{proposition}

\section{Definitions, examples and properties of plt blow-ups}
\begin{definition}
\label{def-plt-blowup}
Let $X$ be a normal algebraic variety and let $f\colon Y\to X$
be a blow-up such that
the exceptional locus of $f$ contains only one irreducible divisor, say $S$.
 Then $f\colon (Y,S)\to X$ is called a \textit{
plt blow-up} of $X$ if $(Y,S)$ is plt and $-(K_Y+S)$ is $f$-ample.
\end{definition}
This definition can be easily modified for blow-ups of
analytic spaces in a usual way (recall that we are working in
the algebraic category).

\begin{remark} (i)
Note that if $X$ is ${\QQ}$-factorial,
then $\rho(Y/X)=1$ and $Y$ is also $\QQ$-factorial.
Indeed, let $D\ne S$ be any effective Weil divisor
on $Y$ intersecting $S$. Then we have
$f^*(f(D))\mathbin{\sim_{\scriptscriptstyle\QQ}}D+\alpha S$,
where $\alpha>0$. Thus
$D\mathbin{\sim_{\scriptscriptstyle\QQ}} -\alpha S$ over $X$.
Since $\mt{Weil}(Y)/\mt{Weil}(X)$ is generated by $S$,
  the variety $Y$ is $\QQ$-factorial.
\par (ii)
If $f\colon (Y,S)\to X$ is a blow-up of a
variety $X$ with only klt
singularities  such that
$\rho (Y/X)=1$ and  $(Y,S)$ is plt, then
$-(K_Y+S)$ is $f$-ample automatically, so $f$
is a plt blow-up.
\par (iii)
Any plt blow-up $f\colon (Y,S)\to X$ is
uniquely defined up to isomorphism over $X$ by the
 discrete valuation corresponding to $S$ (cf. \cite[Proof of 6.2]{Ut}).
\end{remark}

\begin{remark}
By Inversion of Adjunction $S$ is normal and
$(S,\Diff_S(0))$ is klt.
In dimension two $S$ is an
irreducible non-singular curve and since $-(K_S+\Diff_S(0))$
is ample, $S\simeq {\mathbb P}^1$. In dimension three two cases are
possible: $f(S)$ is a point and then $(S,\Diff_S(0))$ is a log
Del Pezzo surface, or $\dim (f(S))=1$ and then $S\to f(S)$ is
generically $\PP^1$-fibration.
\end{remark}

\begin{example}
We describe below the class of all
$(S,\Diff_S(0))$'s for two-dimensional Du Val
singularities (see e.~g. \cite{Pro}). In this case $S\simeq{\mathbb P}^1$
and $\Diff_S(0)=\sum_{i=1}^r (1-1/m_i)P_i$, where $%
P_1,\dots,P_r\in{\mathbb P}^1$ are distinct points, $r\ge 0$ and $m_i\in{%
\NN}$. Since $-(K_S+\Diff_S(0))$ is ample, $\sum
(1-1/m_i)<2$. It is easy to show that for $(m_1,\dots,m_r)$ one has one of
the following possibilities.

\begin{itemize}
\item $r=0$, $1$ or $2$. In this case $(X\ni P)$ is a cyclic quotient
singularity (type $A_n$) and $f$ is toric;
\item $(m_1,\dots ,m_r)=(2,2,m)$, $m$ is arbitrary (type $D_n$);
\item $(m_1,\dots ,m_r)=(2,3,3)$ (type $E_6$);
\item $(m_1,\dots ,m_r)=(2,3,4)$ (type $E_7$);
\item $(m_1,\dots ,m_r)=(2,3,5)$ (type $E_8$).
\end{itemize}

In cases $D_n$ and $E_n$ a plt blow-up is unique
(see \ref{2-unique}).
Note that in two-dimensional case all
(one-dimensional) log Fano varieties with standard boundaries can be
realized as exceptional divisors of some plt blow-up.
Similar description can be obtained for two-dimensional
lt singularities.
\end{example}

We are interested in the three-dimensional case.
\begin{problem}
\label{problem}
Describe the class of all log Del Pezzo surfaces
(resp. generically $\PP^1$-fibrations) which
can be exceptional divisors of some plt blow-up of a
terminal (resp. canonical, $\epsilon$-lt) singularity.
\end{problem}

\begin{example}
\label{w-blow}
Let $(X\ni P)$ be a cyclic quotient singularity
${\CC}^3/{\ZZ}_m$, where $\ZZ_m$ acts on
$\CC^3$ with weights $(a_1,a_2,a_3)$,
 $(a_1,a_2,a_3)=1$.
Then any weighted blow-up $f\colon Y\to X$ is a plt blow-up.
The exceptional divisor $S$ is a weighted
projective space $\PP(a_1,a_2,a_3)$, $\Diff_S(0)$ is
supported in the triangle $\{x_1x_2x_3=0\}$ and has
coefficients $1-\frac{1}{(a_i,a_j)}$.
\end{example}

\begin{example}
Let $\widetilde{{\CC}^3}\to {\CC}^3$ be the blow-up of $x$-axis
and let $\widetilde{{\CC}^3}/{\ZZ}_m\to {\CC}^3/{\ZZ}%
_m$ be a quotient by a cyclic group acting on ${\CC}^3$ free in
codimension $1$. Then
$\widetilde{{\CC}^3}/{\ZZ}_m\to {\CC}^3/{\ZZ}_m$ is a plt blow-up.
The exceptional divisor $S$ is a generically $\PP^1$-bundle
with two singular points on $f^{-1}(0)$.
\end{example}

Combining \ref{prodolj} with \ref{podnyal} (i) we get
the following assertion which shows that
the existence of complements can be established
inductively.

\begin{corollary}
Let $f\colon (Y,S)\to X$ be a plt blow-up of a klt singularity.
Assume that $(S,\Diff_S(0)$ is $n$-complementary.
Then $K_X$ is $n$-complementary.
\end{corollary}

The following proposition shows that plt blow-ups are easy to construct.

\begin{proposition}[cf. {\cite[3.1]{Sh2}}, {\cite[17.10]{Ut}},
 {\cite[5.19]{Sh}}]
\label{per}\label{pltblowup}\label{plt}
Let $X$ be a normal ${\QQ}$-factorial variety of dimension
$\le 3$ and let $D$ be a boundary on $X$
such that $K_X+D$ is lc, but is not plt. Assume that $D\ne 0$
and $X$ has at
worse klt singularities. Then there exists a plt blow-up
$f\colon (Y,S)\to X$ such
that

\begin{enumerate}
\renewcommand\labelenumi{\textrm{(\roman{enumi})}}
\item $K_Y+S+D_Y=f^{*}(K_X+D)$ is lc;
\item $K_Y+S+(1-\varepsilon)D_Y$ is plt and anti-ample over $X$
for any $\varepsilon>0$;
\item
$Y$ is ${\QQ}$-factorial and $\rho (Y/X)=1$.
\end{enumerate}
\end{proposition}

This blow-up $f\colon (Y,S)\to X$ is called a \textit{plt blow-up} of $(X,D)$.
Note that the proposition holds in arbitrary dimension
modulo the Log Minimal Model Program.

\begin{proof}
By Theorem~\ref{min-lt} there exists a blow-up
 $h\colon W\to X$ with the reduced exceptional divisor
$E=\sum E_i$ such that
$W$ is $\QQ$-factorial and $K_W+D_W+E=h^*(K_X+D)$ is dlt.
Note that $\rho(W/X)$ is the number of components of $E$,
because $X$ is $\QQ$-factorial.
Moreover, since $K_X+D$ is not plt, we have $E\ne 0$.

Now we apply the $(K_W+E)$-Minimal Model Program over $X$.
The divisor $K_W+E\equiv h^*(K_X+D)-D_W$ cannot be nef over $X$.
Therefore on the last step we get a divisorial extremal contraction
 $f\colon Y\to X$ with $\rho(Y/X)=1$, with a unique exceptional
divisor $S$ and negative with respect to $K_{Y}+S$.
Because $K_W+E+D_W$ is numerically trivial,
the divisor $f^*(K_X+D)=K_{Y}+S+D_{Y}$ is lc on $Y$.
Since $K_Y+S$ is plt, so is $K_Y+E+(1-\varepsilon)D$
  by \cite[2.17.5]{Ut}.
This completes the proof.
\end{proof}

\begin{proposition}
\label{center}
Let $f\colon (Y,S)\to X$ be a plt blow-up of a variety
$X$ with only klt singularities. Then there exists a
boundary $D$ on $X$ such
that $(S,D)$ is lc, $a(S,D)=-1$ and $S$ is the only divisor
over $X$ with discrepancy
$a(\phantom{S},D)\le -1$.
\end{proposition}

\begin{proof}
For sufficiently big and divisible $n$ the divisor
$-n(K_Y+S)$ is very ample over $X$.
Take a general member $L\in |-n(K_Y+S)+f^*(H)|$,
where the divisor $H$ is ample enough on $X$.
Put $B:=(1/n)L$,
$D:=f(B)$. Then $(Y,S+B)$ is plt and since
$f^*(K_X+D)=K_Y+S+B$, the pair $(X,D)$ is lc
(see \cite[3.10]{KoP}).
\end{proof}

\begin{corollary}
Notation as in Proposition~\ref{center}. Then

\begin{enumerate}
\item
$f(S)$ is the minimal center of $LCS(X,D)$;
\item
$f(S)$ is normal;
\item
if $\codim (f(S))=2$, then there exists a
boundary $D_{f(S)}$ on $f(S)$ such that $(f(S),D_{f(S)})$ is klt.
\end{enumerate}
\end{corollary}

\begin{proof}
Since $S$ is the only divisor over $X$ with discrepancy $a(S,f(B))=-1$,
$f(S)$ is the minimal center of lc singularities of $(X,D)$.
It is normal by \cite{K1}, (iii) follows by \cite{K2}.
\end{proof}

The following is converse to Proposition~ \ref{center}.

\begin{proposition}
\label{2equiv}
Let $X$ be a normal $\QQ$-factorial
variety of dimension $\le 3$ and let
$W\subset X$ be an irreducible subvariety.
Assume there exists a boundary $D$
such that $(X,D)$ is lc and the minimal center of $LCS(X,D)$ is $W$. Then
there exists a plt blow up with center along $W$.
\end{proposition}

\begin{proof}
Apply Proposition~\ref{pltblowup} to $(X,D)$.
\end{proof}

\begin{corollary}
\label{curve-point}
Let $X$ be a normal three-dimensional variety and let $f\colon (Y\supset
S)\to X$ be a plt blow-up with center along a curve $W=f(S)$. Then for any
point $P\in W$ there exists a plt blow-up with center at $P$.
\end{corollary}

\section{The terminal case}
\begin{proposition}
\label{m}
Let $(X\ni P)$ be a three-dimensional
 klt singularity and let
$f\colon (Y,S)\to X$ be a plt blow-up. Assume either
\begin{enumerate}
\item
$\dim f(S)=1$, or
\item
there
exists a boundary $D$ on $X$ such that $\down{D}\ne \emptyset$ near
$P$ and $(X,D)$ is plt.
\end{enumerate}
Then $K_S+\Diff_S(0)$ has a regular non-klt complement
(near the fiber $f^{-1}(P)$ in the case $\dim f(S)=1$).
\end{proposition}

\begin{proof}
The assertion (i) follows from Proposition~\ref{P-bundle}.
To prove (ii) we assume first that
$X$ is $\QQ$-factorial.
Write $f^*(K_X+D)=K_Y+D_Y+aS$, where $D_Y$ is the proper transform
of $D$ on $Y$ and $a<1$.
By Connectedness Lemma~\cite[5.7]{Sh}, \cite[17.4]{Ut}
the set of lc singularities $LCS(Y, D_Y+S)$ is connected.
It is clear that $LCS(Y, D_Y+S)$ contains $S$ and $\down{D_Y}$.
Whence $K_Y+D_Y+S$
is not plt (see \cite[3.6]{Sh}).
Since $K_Y+S$ is plt, $K_Y+bD_Y+S$ is lc but is not
plt for
some $0<b\le 1$. Therefore $K_S+\Diff_S(bD_Y)$ is lc but not
klt (see Theorem~\ref{Inv-Adj} and
\cite[17.7]{Ut}). On the other hand $-(K_S+\Diff_S(bD_Y))$ is ample.
The assertion follows by Proposition~\ref{DelPezzo}.
If $X$ is not $\QQ$-factorial, then we consider
a $\QQ$-factorialization $h\colon Y'\to Y$ (see \cite[6.11.1]{Ut}) and let
$f'\colon (Y',S') \to X$ be the composition map.
The same arguments show that $K_{S'}+\Diff_{S'}(0)$
has a regular non-klt complement $K_{S'}+\Diff_{S'}(0)^+$
(we have to remark that in this case $-(K_{S'}+\Diff_{S'}(0))$
is nef and big because $h$ is a small morphism).
Finally, to conclude the proof we can
push $K_{S'}+\Diff_{S'}(0)^+$ down
on $S$ by (i) of Proposition~\ref{podnyal}.
\end{proof}

\begin{corollary}
In notation of Proposition~\ref{m} $K_Y+S$ has a
regular non-plt complement.
\end{corollary}

The following is a necessary condition
for Problem~\ref{problem} in the terminal case.

\begin{corollary}
\label{m1}
Let $(X\ni P)$ be a three-dimensional
 terminal singularity and
let $f\colon (Y,S)\to X$ be a plt blow-up. Then
$K_S+\Diff_S(0)$ has a regular non-plt complement.
\end{corollary}

\begin{proof}
Consider $f\colon (Y\supset S)\to (X\ni P)$ as a morphism of
analytic spaces.
By \cite[6.4.B]{YPG} the general element $D\in |-K_X|$
has only a Du Val singularity at $P$. Inversion
of Adjunction gives us that $(X,D)$ is plt, so
to complete the proof we can apply
the analytic analog of Proposition~\ref{m}.
\end{proof}
Note that \ref{m1} has no generalizations for
canonical singularities (see \cite{MP}, \cite{MP1} and
Example~\ref{ex-2}).
\par
Using the same arguments as in the proof of
Proposition~\ref{m} one can show the following

\begin{proposition}
\label{mmm}
Let $(X\ni P)$ be a three-dimensional
 klt singularity, let
$f\colon (Y,S)\to X$ be a plt blow-up, and let
$K_X+D$ be an $n$-complement near $P$.
Then one of the following holds
\begin{enumerate}
\item
$a(S,D)=-1$ and $K_Y+S+D_Y:=f^*(K_X+D)$ is an $n$-complement of
$K_Y+S$;
\item
$a(S,D)>-1$ and then there exists a regular complement of
$K_Y+S$ which is not plt.
\end{enumerate}
\end{proposition}

\begin{corollary}
\label{mmm1}
Let $(X\ni P)$ be a three-dimensional klt singularity, let
$f\colon (Y,S)\to X$ be a plt blow-up.
Assume that $K_X$ has a regular complement $K_X+D$.
Then $K_Y+S$ also has a regular complement $K_Y+S+B$
such that $K_Y+S+B=f^*(K_X+D)$ or $K_Y+S+B$ is not plt.
\end{corollary}

\begin{example}
\label{ex-1}
Let $(X\ni P)$ be the terminal $cE_8$-singularity given by
the equation $x_1^2+x_2^3+x_3^5+x_4^5=0$ and let
$f\colon (Y,S)\to X$ be the weighted blow-up
with weights $(15,10,6,6)$. Then $S$ is isomorphic to the cone over a
rational normal curve of degree $5$ and
$\Diff_S(0)=(1/2)H_1+(2/3)H_2$, where $H_1$, $H_2$ are
hyperplane sections in a
general position.
Since $(S,(1/2)H_1+(2/3)H_2)$ is klt, $f\colon (Y,S)\to X$
is a plt blow-up.
It is easy to check that $-K_S=(7/5)H$.
Since $\Diff_S(0)^+\ge\Diff_S(0)$,
$K_S+\Diff_S(0)$ is not $1$ or $2$-complementary.
On the other hand there exists a $3$-complement
$K_S+(2/3)H_1+(2/3)H_2+1/3L$, where $L$ is a line on $S$.
Note however that this complement is klt.
The only regular non-klt $n$-complement is
$K_S+(1/2)H_1+(2/3)H_2+L_1+(1/6)L_2$, where
$L_1$, $L_2$ are lines on $S$. Here  $n=6$.
\end{example}

Unfortunately the proof of \cite[6.4.B]{YPG} uses the classification of
terminal singularities and has no generalizations in higher dimensions.
To avoid this difficulty one can apply the method proposed by
Shokurov in \cite[7.8]{Sh1}. This method uses the classification of
surface complements and computations of discrepancies.
Remark also that the necessary condition of \ref{m1} is
not sufficient. Below we give another necessary condition
for Problem~\ref{problem}.

\begin{proposition}
\label{cond2}
Let $f\colon (Y,S)\to X$ be a plt blow-up of
a terminal singularity of index $r$
and let $L$ be a curve on $S$ such that
$L\cap\Sing(S)=L\cap\Sing(\Supp(\Diff_S(0)))=\emptyset$.
Write $\Diff_S(0)=\sum (1-1/m_i)\Delta_i$, where $m_i\in\NN$
(see \ref{coeff}) and denote $m:=\mt{lcm}\{ m_i\}$.
Then
\begin{equation}
\label{t-eq}
-m(K_S+\Diff_S(0))\cdot L\ge 1+1/r.
\end{equation}
\end{proposition}
Note that $-m(K_S+\Diff_S(0))\cdot L$ can be computed on $S$,
without using any information about whole $Y$.
\begin{proof}
We can write
$$
K_Y+S=f^*K_X+(1+a)S,\quad\text{where}\ a=a(S,0,X).
$$
Then along $L$ divisors
$m(K_Y+S)$ and $mS$ are Cartier
by the following lemma (cf. \ref{coeff}).

\begin{lemma}[{\cite{Pro2}}]
\label{odna}
Let $(Y,S)$ be a plt pair with a reduced boundary $S\ne 0$.
Assume that near some point $Q$
the surface $S$ is non-singular and
$\Diff_S(0)=(1-1/m)\Delta_1$, where $\Delta_1$
is a non-singular curve.
Moreover, assume that $Y$ is nonsingular outside $S$.
Then near this point
there is an analytic isomorphism
$(Y,S)\simeq (\CC^3,\{x_1=0\})/\ZZ_{m}(0,1,a)$, where $(a,m)=1$.
\end{lemma}
\begin{proof}
Let $\varphi\colon Y'\to Y$ be a cyclic cover such that
$S':=\varphi^{-1}(S)$ is a Cartier divisor (see \cite[3.6]{YPG}).
Then $(Y',S')$ is plt by \cite[\S 2]{Sh}, \cite[16.3]{Ut}.
Therefore the singularity $Y'$ at $Q':=\varphi^{-1}(Q)$
is terminal. Consider its canonical cyclic cover
$\psi\colon Y''\to Y'$ and let $S'':=\psi^{-1}(S')$.
The finite morphism of normal surfaces $S''\to S$
is \'etale outside $\Supp(\Sing(Y))=\Supp(\Diff_S(0))$.
Since the fundamental group $\pi_1(S\setminus\Supp(\Diff_S(0))$
is isomorphic to $\ZZ$ (because
$(S,\Supp(\Diff_S(0)))\simeq (\CC^2,\CC^1)$),
we have that $\psi\colon Y''\to Y'$ is a cyclic cover.
Moreover, $S''\to S$ has a smooth ramification divisor.
Hence $S''$ is non-singular and so is $Y''$.
The rest is obvious.
\end{proof}

Thus we have
$$
-m(K_S+\Diff_S(0))\cdot L=
-m(1+a)S\cdot L,
$$
where $-m(K_S+\Diff_S(0))\cdot L=-m(K_Y+S)\cdot L$
and $-mS\cdot L$ are positive integers.
Since $(X\ni P)$ is terminal of index $r$, we have
$a\ge 1/r$. This gives us the desied inequality (\ref{t-eq}).
\end{proof}

More delicate method was proposed by Shokurov in \cite[7.8]{Sh1}.
It uses computations of all discrepancies of $(X\ni P)$.

\begin{example}
Assume that $S$ is non-singular and $\Diff_S(0)=\emptyset$
(so $S$ is a usual Del Pezzo surface).
Then any curve on $S$ satisfies conditions of
 \ref{cond2}. By (\ref{t-eq})  $-K_S\cdot L>1$, so we have
only two cases: $S=\PP^2$ or $S=\PP^1\times\PP^1$.
Note that in these cases $Y$ is non-singular and
$Y\to X$ is  Mori's extremal contraction.
Similar, but more complicated, description can be done in
the case when $S$ has only Du Val singularities.
In particular, then we have $K_S^2\ne 1$ and $7$.
\end{example}

Corollary~\ref{m1} and Proposition~\ref{cond2} show that
log Del Pezzo surfaces which are exceptional
divisors of plt blow-ups of a terminal
singularities form a very restricted class.
Regular complements of surfaces were classified in \cite[\S 6]{Sh1}. They
are divided into types $\AA_m^n$, ${\mathbb E}1_m^n$,
${\mathbb D}_m^n$, ${\mathbb E}2_m^n$, ${\mathbb E}3_m^n$,
${\mathbb E}4_m^n$ and
${\mathbb E}6_m^n$.
It is a very interesting problem to establish a correspondence
between this classification and classification of
terminal singularities \cite{YPG}. For example, the complements
$(S,\Diff_S(0)^+)$
from \ref{ex-1} are of types ${\mathbb E}3_0^0$
and ${\mathbb E}6_0^1$.

\section{On uniqueness of plt blow-ups}
In this section we discuss the uniqueness property of plt blow-ups.

\begin{definition}
\label{def-exc}\label{def-weakly-exc}
Let $(X\ni P)$ be a klt singularity. It said to be \textit{weekly
exceptional} if there exists only one plt blow-up
(up to birational equivalence).
A lc pair $(X,D)$ is said to be \textit{exceptional}
if there exists at most one divisor $E$ over $X$ with discrepancy
$a(E,D)=-1$.
A normal klt singularity $(X\ni P)$ is said to be \textit{exceptional}
if $(X,D)$ is exceptional for any boundary $D$
whenever $K_X+D$ is log canonical (see {\cite[1.5]{Sh1}}).
\end{definition}

\begin{proposition}[\cite{MP1}]
\label{convex}
Let $(X\ni P)$ be an exceptional $\QQ$-factorial
singularity. Then the
divisor $S$ from Definition \ref{def-exc} with
discrepancy $a(S,D)=-1$ is
birationally unique. Therefore if a singularity is
exceptional, then it is also weakly exceptional.
\end{proposition}

Below we give criterions for a klt singularity to be
(weakly) exceptional in terms of plt blow-ups.

\begin{theorem}
\label{weakly-exc}
Let $(X\ni P)$ be a $\QQ$-factorial klt singularity
of dimension $\le 3$ and let
$f\colon (Y,S)\to X$
be a plt blow-up of $P$. Then the
following are equivalent

\begin{enumerate}
\item
$(X\ni P)$ is not weakly exceptional;

\item
there is an effective $\QQ$-divisor $\Theta\ge\Diff_S(0)$
such that $-(K_S+\Theta)$ is ample and $(S,\Theta)$ is not lc;

\item
there is an effective $\QQ$-divisor $\Theta\ge\Diff_S(0)$
such that $-(K_S+\Theta)$ is ample and $(S,\Theta)$ 
is not klt\footnote{The author would like to thank
Dr. O. Fujino, who pointed out some mistakes in this theorem}.
\end{enumerate}
\end{theorem}

\begin{proof}
Note that $\rho(Y/X)=1$ and $Y$ is $\QQ$-factorial
(because so is $X$).
First we show the implication (i) $\Longrightarrow$ (ii).
Assume that $(X\ni P)$ is not weakly exceptional.
By Proposition~\ref{center} there exists a boundary $D$
such that $(X,D)$ is lc and $a(S,D)>-1$, $a(E,D)=-1$
for some $E\ne S$. Let us write
$K_Y+\alpha S+D_Y=f^*(K_X+D)$, where $D_Y$ is the proper transform
of $D$ and $\alpha<1$.
This divisor is lc but is not plt.
Then $-(K_Y+S+D_Y)$ is ample over $X$ and $(Y,S+D_Y)$
is not lc. Take $\Theta=\Diff_S(D_Y)$.
\par
The implication (ii) $\Longrightarrow$ (iii) is clear.
To show
(iii) $\Longrightarrow$ (ii) we replace $\Theta$ on
$\QQ$-Cartier divisor $\Theta+\epsilon(\Theta-\Diff_S(0))$
for some small $\epsilon \ge 1$.
\par
Finally, we prove (ii) $\Longrightarrow$ (i).
By \ref{curve-point} we may assume that $f(S)$ is a point.
For sufficiently large and divisible $n$
the integer Cartier divisor
$-n(K_S+\Theta)$ is a very ample.
Take $H\in |-n(K_S+\Theta)|$ and denote $\Xi:=\Theta+\frac{1}{n}H$.
Then $K_S+\Xi$ is numerically trivial and is not lc.
Moreover, we may assume that $n(K_S+\Xi)\sim 0$.
From the exact sequence
\begin{multline*}
0\longrightarrow\OOO_Y(K_Y-(n+1)(K_Y+S))
\longrightarrow\OOO_Y(-n(K_Y+S))\\
\longrightarrow
\OOO_S(-n(K_S+\Diff_S(0)))
\longrightarrow 0
\end{multline*}
and the Kawamata-Viehweg vanishing
$H^1(Y,\OOO_Y(K_Y-(n+1)(K_Y+S)))=0$ we get that
the map
$$
H^0(Y,\OOO_Y(-n(K_Y+S)))
\longrightarrow
H^0(S,\OOO_S(-n(K_S+\Diff_S(0))))
$$
is surjective. Thus
there exists an integer effective divisor
$L\in |-(K_Y+S)|$ such that $L|_S=n(\Xi-\Diff_S(0))$.
Consider $\QQ$-Cartier divisor $B:=\frac{1}{n}L$. Then
$n(K_Y+S+B)\sim 0$ and
$(K_Y+S+B)|_S=K_S+\Xi$. By \cite[17.7]{Ut}
$K_Y+S+B$ is not lc. We can take $B'=\alpha B$, $\alpha<1$ such that
 $K_Y+S+B'$ is lc but is not plt. It is easy to see that
\begin{equation}
\label{eq1}
-(K_Y+S+B')\equiv (1-\alpha)B\equiv -(1-\alpha)(K_Y+S)
\end{equation}
is $f$-ample.
Further, we take a minimal log terminal blow-up (see \ref{min-lt})
$g\colon W\to Y$. We have
$K_W+S_W+E+B'_W=g^*(K_Y+S+B')$, where $S_W$ and $B'_W$ are
proper transforms of $S$ and $B'$, respectively, and
$E\ne \emptyset$ is the reduced exceptional divisor.
One can see by (\ref{eq1})
\begin{equation}
\label{eq2}
K_W+S_W+E+B'_W\equiv -(1-\alpha)(g^*(B).
\end{equation}

Apply the Log Minimal Model Program to $(W,S_W+E+B'_W)$
over $X$.
On each step by (\ref{eq2}) $K_W+S_W+E+B'_W$
cannot be nef over $X$. Let $V\to X$
be the last step (i.~e.  $V\to X$
contracts an irreducible exceptional divisor, say $R$).
By construction  $K_V+R+B'_V$ is plt.
Since $X$ is $\QQ$-factorial,
$\rho(V/X)=1$ and $(V,R)\to X$ is a plt blow-up.
It is sufficient to show that $R$ is the proper
transform of some of components of $E$
(otherwise $(X\ni P)$ is not weakly exceptional). Indeed,
in the opposite case $(V,R)\to X$ is birationally
equivalent to $(Y,S)\to X$. Hence they are isomorphic
because $\rho(Y/X)=\rho(V/X)=1$.
But then $K_Y+S+B'$ is plt.
This contradiction concludes the proof.
\end{proof}

\begin{corollary}
\label{weakly-cor}
Notation as in Theorem~\ref{weakly-exc}.
If a three-dimensional singularity  $(X\ni P)$
is not weakly exceptional, then $K_S+\Diff_S(0)$ has a regular
non-klt complement $K_S+\Diff_S(0)^+$. Moreover,
one can choose $\Diff_S(0)^+$ in such a way that
$a(E, \Diff_S(0)^+)=-1$, $a(E,\Theta)\le -1$ for some
(not necessary exceptional divisor $E$).
\end{corollary}
\begin{proof}
Apply Proposition~\ref{DelPezzo} or \ref{P-bundle}.
\end{proof}

\begin{corollary}
Notation as in \ref{weakly-exc} and \ref{weakly-cor}.
If a three-dimensional singularity
$(X\ni P)$ is not weakly exceptional, then $K_X$ has a regular
complement. Moreover, $K_Y+S$ has a regular non-plt complement.
\end{corollary}

\begin{corollary}
\label{weakly-exc-red}
Let $(X\ni P)$ be a klt singularity and let
 $D$ be a boundary on
$X$ such that $\down{D}\ne \emptyset$ near $P$.
If $(X,D)$ is plt and $\dim(X)\le 3$, then $(X\ni P)$ is not
weakly exceptional.
\end{corollary}
\begin{proof}
Let $f\colon (Y,S)\to X$ be a plt blow-up.
Write $K_Y+D_Y+\alpha S=f^*(K_X+D)$, where $\alpha<1$ and $D_Y$ is
the proper transform of $D$.
Then $K_Y+D_Y+S$ is anti-ample over $X$ and it is not plt.
Therefore in Theorem~\ref{weakly-exc} we can take
$\Theta=\Diff_S(D_Y)$.
\end{proof}

\begin{example}
\label{2-unique}
Let $(X\ni P)$ be a two-dimensional Du Val
singularity. By \ref{weakly-exc-red} and by
\cite[\S 9]{K} it is not weakly exceptional
only if it is of type $A_n$. On the other hand
singularities of type $A_n$ are not weakly
exceptional (cf. \ref{w-blow}). In cases $D_n$,
$E_6$, $E_7$ and $E_8$ the unique plt blow-up
is the blow-up of the "central" vertex of
the dual graph.
\end{example}
Similar to Proposition\ref{m} one can prove
\begin{corollary}
\label{term-non-un}
Let $(X\ni P)$ be a three-dimensional terminal singularity.
Then it is not weakly exceptional.
\end{corollary}

Now we prove a criterion of the exceptionality.
In contrast to Theorem~\ref{weakly-exc} it does not
depend on the Minimal Model Program.

\begin{theorem}
\label{exc}
Let $(X\ni P)$ be a klt singularity and let
$f\colon (Y,S)\to X$
be a plt blow-up of $P$. Then the
following are equivalent

\begin{enumerate}
\item $(X\ni P)$ is non-exceptional;
\item there is an effective $\QQ$-divisor
$\Theta\ge\Diff_S(0)$ such that $n(K_S+\Theta)\sim 0$
for some $n\in\NN$ and   $K_S+\Theta$ is not klt;
\item there is an effective $\QQ$-divisor 
$\Xi\ge\Diff_S(0)$ such that
$-(K_S+\Xi)$ is nef (in particular $\QQ$-Cartier)
and $(S,\Xi)$ is not klt;
\item
(in dimension $3$ only)
there is a regular (i.~e. $1$, $2$, $3$, $4$ or $6$)
complement of $K_S+\Diff_S(0)$ which is not klt;
\item
(in dimension $3$ only)
there is a regular complement $K_X+D$ of $K_X$ which is not
exceptional and such that $a(S,D)=-1$ (and then the set of
divisors with discrepancy $-1$ must be infinite).
\end{enumerate}
\end{theorem}

\begin{proof}
Note that $S$ and $K_Y$ are $\QQ$-Cartier,
because $(X\ni P)$ is klt.
First we prove (ii) $\Rightarrow$ (i).
Assume that $(X\ni P)$ is exceptional.
Let $\Theta$ be a boundary
such as in (ii). Then $\Theta\ge \Diff_S(0)$ and $K_S+\Theta\equiv 0$.
Take $n\in \NN$ such that
$n(\Theta-\Diff_S(0))$
and $n(K_Y+S)$ are integer Cartier divisors. It is clear that
$n(\Theta-\Diff_S(0))\in |-n(K_S+\Diff_S(0))|$.
We have the exact sequence
\begin{multline*}
0\longrightarrow\OOO_Y(K_Y-(n+1)(K_Y+S))
\longrightarrow\OOO_Y(-n(K_Y+S))\\
\longrightarrow
\OOO_S(-n(K_S+\Diff_S(0)))
\longrightarrow 0
\end{multline*}
By the Kawamata-Viehweg vanishing
$H^1(Y,\OOO_Y(K_Y-(n+1)(K_Y+S)))=0$, because $-(K_Y+S)$
is ample. So there exists an integer
effective Weil divisor $L$ on $Y$ such that
$(K_Y+S+\frac{1}{n}L)_S=K_S+\Theta$.
Denote $B:=\frac{1}{n}L$. If $K_S+\Theta$ is not klt, then
by Inversion of Adjunction $K_Y+S+B$ is not plt.
Put $D:=f(B)$.
Thus we have either $K_X+D$ is lc (and then there
is an exceptional divisor
$E\ne S$ over $X$ with $a(E,D)=a(S,D)=-1$) or
$K_X+\alpha D$ is lc but is not klt
for some $\alpha<1$. Since $a(S,D)=-1$, there
is an exceptional divisor $E\ne S$ over $X$ with $a(E,\alpha D)=-1$.
Both cases gives us a contradiction.
\par
The implications (ii) $\Rightarrow$ (iii) and
(iv) $\Rightarrow$ (ii) are obvious. As for (iii) $\Rightarrow$ (ii)
we can take $\Theta=\Xi+(1/n)F$, where $F\in |-n(K_S+\Xi)|$,
the linear system $|-n(K_S+\Xi)|$ is not empty
for $n\gg 0$ by Base Point Free Theorem \cite[3-1]{KMM}.
\par
Let us show (i) $\Rightarrow$ (ii).
By the assumption $(X\ni P)$ is not exceptional.
Therefore there exists a boundary $D$
such that $(X,D)$ is lc and $a(E,D)=-1$
for some $E\ne S$. Let us write
$K_Y+\alpha S+B=f^*(K_X+D)$, where $\alpha\le 1$.
Then $-(K_Y+S+B)$ is nef over $X$ and $(Y,S+B)$
is not plt. Take $\Theta=\Diff_S(B)$.
\par
To see (v) $\Rightarrow$ (iv) we can take the regular complement
$f^*(K_X+D)|_S$ on $K_S+\Diff_S(0)$ which is not klt by
Inversion of Adjunction~\ref{Inv-Adj}. The converse implication
follows from Propositions~\ref{prodolj} and \ref{podnyal} (i).
\par
Finally, we prove (ii) $\Rightarrow$ (iv).
We may assume that $f(S)$ is a point (otherwise there exists
a regular non-klt complement by Proposition~\ref{P-bundle}).
Thus $S$ is projective.
Denote $\Diff_S(0)$ by $\Delta=\sum\delta_i\Delta_i$.
Assume that there exists a boundary $\Theta\ge\Delta$ on $S$ such that
$K_S+\Theta$ is numerically trivial and is not klt.
The implication (ii) $\Rightarrow$ (iv)
is an immediate consequence of
Corollary~\ref{DelPezzo2} below which is a
 generalization of \ref{DelPezzo}.
\end{proof}

\begin{theorem}[{\cite[2.3]{Sh1}}]
\label{DelPezzo1}
Let $S$ be a projective surface and let $\Theta$
be a boundary on it with standard coefficients
such that $-(K_S+\Theta)$ is nef. Assume that $(S,\Theta)$
is lc but is not klt and assume that the Mori cone
$\NE(S)$ is polyhedral and has contractible faces.
Then $K_S+\Theta$ has a regular complement.
\end{theorem}

\begin{corollary}
\label{DelPezzo2}
Let $S$ be a projective surface and let $\Theta\ge\Delta$
 be a boundaries on it.
Assume that
\begin{enumerate}
\item
$\Delta$ has only standard coefficients;
\item
$(S,\Delta)$ is klt;
\item
$-(K_S+\Delta)$ is ample;
\item
$-(K_S+\Theta)$ is nef;
\item
$(S,\Delta)$ is lc but is not klt.
\end{enumerate}
Then $K_S+\Theta$ has
a regular complement which is not klt.
\end{corollary}
\begin{proof}
Note that $\NE(S)$ is polyhedral and has contractible faces
by Cone Theorem \cite[4-2-1]{KMM}.
First we consider the case when $\down{\Theta}\ne 0$.
Let $\Delta_1$ be a component of $\down{\Theta}$.
Replace $\delta_1$ with $1$:
$\Delta':=\Delta_1+\sum_{i\ne 1}\delta_i\Delta_i$, so
 $\Delta\le\Delta'\le\Theta$.
If $-(K_S+\Delta')$ is nef, then we can apply \ref{DelPezzo1}
(because $\down{\Delta'}\ne 0$ and $\Delta'$ has
only standard coefficients).
Further, we assume that $-(K_S+\Delta')$ is not nef.
Then there exists a $(K_S+\Delta)$-extremal ray $R$ such that
$(K_S+\Delta')\cdot R\ge 0$. If it is birational, then we
contract it. Since $K_S+\Theta$ is trivial on $R$, it preserves
the lc property of $K_S+\Theta$ and $K_S+\Delta'$. Since
$\Delta'$ has only standard coefficient, we can pull-back
complements of $K_S+\Delta'$
(but here we are looking for regular complements of
$K_S+\Delta'$, see Lemma~\ref{podnyal}).
Note also that $(\Delta'-\Delta)\cdot R>0$.
Therefore the curve $\Delta_1$ is not contracted.
Continue the process we get a
non-birational extremal ray $R$ on $S$
such that $(K_S+\Delta')\cdot R> 0$.
But on another hand,
$(K_S+\Delta')\cdot R\le (K_S+\Theta)\cdot R=0$, a contradiction.
\par
Let now $\down{\Theta}=0$, then at some point, say $P$, the divisor
$K_S+\Theta$ is lc, but is not plt.
In this point $K_S+\Delta$ is klt.
Similar to Proposition~\ref{per} we can construct a blow-up
$f\colon\widetilde{S}\to S$ with an irreducible
exceptional divisor $E$ over $P$
 such that $a(E,\Theta)=-1$, the crepant pull back
$K_{\widetilde{S}}+\widetilde{\Theta}+E=f^*(K_S+\Theta)$
is lc and $K_{\widetilde{S}}+\widetilde{\Theta}'+E$
 is plt for any
$\widetilde{\Theta}'=\widetilde{\Delta}+
\epsilon(\widetilde{\Theta}-\widetilde{\Delta})$, $\epsilon<1$.
Recall that $\widetilde{\Theta}$ is the proper transform of $\Theta$.
Note also that $\rho(\widetilde{S}/S)=1$.
$$
K_{\widetilde{S}}+\widetilde{\Delta}+\alpha E=f^*(K_S+\Delta),
$$
where $\alpha<1$.
Assume that
there exists a curve $C$ such that
$(K_{\widetilde{S}}+\widetilde{\Delta}+E)\cdot C>0$. Then
$(\widetilde{\Theta}-\widetilde{\Delta})\cdot C<0$.
Therefore $C$ is a component of
$(\widetilde{\Theta}-\widetilde{\Delta})$ and $C^2<0$.
Further, $C\cdot E>0$. Hence $C\ne E$ and we
can chose $\widetilde{\Theta}'<\widetilde{\Theta}$
so that $(K_{\widetilde{S}}+\widetilde{\Theta}')\cdot C<0$.
Therefore the $C$ is an
$(K_{\widetilde{S}}+\widetilde{\Theta}')$-extremal curve
and its contraction preserves the lc property of
$K_{\widetilde{S}}+\widetilde{\Theta}+E$.
Again we can pull-back complements of
$K_{\widetilde{S}}+\widetilde{\Delta}+E$
(see Lemma~\ref{podnyal}).
Continue the process we get the situation when
$-(K_{\widetilde{S}}+\widetilde{\Delta}+E)$
is nef.
All the steps preserve the nef and big property of
$-(K_{\widetilde{S}}+\widetilde{\Delta}+\alpha E)$.
Now we can apply Shokurov's Theorem~\ref{DelPezzo1}.
It is sufficient to prove
that the cone $\NE(\widetilde{S})$ is polyhedral and its faces
are contractractible.
This an easy consequence of the following easy lemma
which is a corollary of Cone Theorem~\cite[4-2-1]{KMM}.

\begin{lemma}[{\cite[2.5]{Sh1}}]
\label{COne}
Let $S$ be a projective surface and let $D$ be a boundary on $S$
such that $(S, D)$ is klt and $-(K_S+D)$ is nef and big
(in particular, $(S, D)$ is a weakly log Del Pezzo). Then the Mori
cone $\NE(S)$ is polyhedral and its faces are contractible.
\end{lemma}

\begin{proof}
By Base Point Free Theorem~\cite[3-1-1]{KMM} the linear system
$|-n(K_S+D)|$ for $n\gg 0$ gives a birational contraction $S\to S'$.
Let $\sum E_i$ be the (reduced) exceptional divisor.
For some small positive $\epsilon_i$'s the
$\QQ$-divisor $-\sum \epsilon_i E_i$, is ample over $S'$.
The assertion follows by Cone Theorem~\cite[4-2-1]{KMM}
applied to $(S,D+ \sum \epsilon_i E_i)$.
\end{proof}

Further, if $\alpha\ge 0$, then we apply Lemma to
$(\widetilde{S},\widetilde{\Delta}+\alpha E)$.
If $\alpha<0$, then by the monotonicity
$-(K_{\widetilde{S}}+\widetilde{\Delta})$
is nef and big. Again apply Lemma~\ref{COne}.
This concludes the proof of Corollary~\ref{DelPezzo2}.
\end{proof}

\begin{corollary}[\cite{Sh1}]
Let $(X\ni P)$ be a three-dimensional
 klt non-exceptional singularity.
Then it is $1$, $2$, $3$, $4$ or $6$ complementary.
\end{corollary}
\begin{corollary}
Let $(X\ni P)$ be a three-dimensional
 klt singularity. Then it is non-exceptional
if and only if there exists a regular complement
$K_X+D$ such that the set $\{E\}$ of divisors with
discrepancy $a(E,D)=-1$ is infinite.
\end{corollary}

\begin{corollary}[{\cite{MP}}]
\label{exc-red}
Let $(X\ni P)$ be a klt singularity and let $D$ be a boundary on
$X$ such that $\down{D}\ne \emptyset$ near $P$.
If $(X,D)$ is lc, then $(X\ni P)$ is not exceptional,
\end{corollary}

V.~V.~Shokurov pointed out to the author that
one can replace the condition $\down{D}\ne 0$ in
Corollaries~\ref{weakly-exc-red} and \ref{exc-red}
on $\max\{d_i\}\ge 41/42$ (see \cite{Ko1}).

\begin{example}[{\cite[5.2.3]{Sh}}]
Let $(X\ni P)$ be a two-dimensional Du Val
singularity. By \ref{exc-red} and by
\cite[\S 9]{K} it is exceptional
if only if it is of type $E_6$, $E_7$ or $E_8$.
\end{example}

Theorems~\ref{weakly-exc} and \ref{exc}
allow us to construct many examples of three-dimensional
exceptional and weakly exceptional singularities singularities
(cf. \cite{MP}, \cite{MP1}).

\begin{example}
\label{ex-2}
Let $X\subset {\CC}^4$ be the hypersurface singularity
$x_1^3+x_2^4+x_3^4+x_4^4=0$. It is canonical by
\cite{R1}. Consider the weighted blow-up $(Y,S)\to X$
with weights $(4,3,3,3)$. The exceptional divisor $S$ is defined
into $\PP(4,3,3,3)$ by the equation $x_1^3+x_2^4+x_3^4+x_4^4=0$.
It is easy to see that the projection $S\to\PP(3,3,3)$
is an isomorphism, so $S\simeq {\mathbb P}^2$. The variety
$Y$ is singular along the curve $x_1=x_2^4+x_3^4+x_4^4=0$ and
along this curve $(Y\supset S)$ is isomorphic to
$(\CC^3\supset\CC^2)/\ZZ_3$. Therefore $(Y,S)$ is plt and
$(Y,S)\to X$ is plt blow-up. Moreover,
 $\Diff_S(0)=(2/3)C$, where $C$ is
a non-singular curve of degree $4$ on $S=\PP^2$.
Among regular complements there are only the following
\begin{itemize}
\item
$3$-complement $K_{\PP^2}+(2/3)C+(1/3)L$, where $L$ is a line,
\item
$6$-complement
$K_{\PP^2}+(2/3)C+(1/6)Q$, where $Q$ is a conic, and
\item
$6$-complement $K_{\PP^2}+(2/3)C+(1/6)L_1+(1/6)L_2$,
where $L_1$, $L_2$ are lines.
\end{itemize}
All these complements are klt by the simple lemma below.
Therefore the singularity is exceptional by Theorem~\ref{exc}.
\end{example}

\begin{lemma}
Let $D=\sum d_iD_i$ be a boundary on ${\CC}^2$, where $D_i$ are
non-singular curves. If $\sum d_i\le 1$ (resp. $<1$), then
$({\CC}^2,D)$ is canonical (resp. terminal).
\end{lemma}

\begin{example}
Let $X\subset {\CC}^4$ be the canonical hypersurface singularity
$x_1^3+x_2^3+x_3^3+x_4^n=0$, where $(n,3)=1$, $n\ge 7$.
 The hyperplane section $x_4=0$ is lc.
Therefore by Corollary~\ref{exc-red} $(X\ni 0)$
is not exceptional.
Consider the weighted blow-up $(Y,S)\to X$
with weights $(n,n,n,3)$.
As in Example~\ref{ex-2} we have $S\simeq\PP^2$, $(Y,S)$ is plt and
 $\Diff_S(0)=(1-1/n)C$, where $C$ is a smooth cubic curve.
By Theorem~\ref{weakly-exc} and Proposition~\ref{DelPezzo}
 the singularity is weakly exceptional.
\end{example}

\end{document}